\documentclass[titlepage]{article}
\usepackage[T2A]{fontenc}
\usepackage[russian]{babel}
\usepackage{amsfonts,amssymb}
\usepackage{color}

\newtheorem{theorem}{Теорема}
\newtheorem{lemma}{Лемма}

\newtheorem{cor}{Следствие}

\begin{document}

\begin{center}
{\bf Четырехмерная простая алгебра с дробной PI-экспонентой}\\

\vskip 0.1in

Зайцев М.В., Реповш Д.\footnote{Работа поддержана РФФИ, грант 13-01-00234а и
Словенским Исследовательским Агенством, грант P1-0292-0101}
\end{center}

{\small
В работе изучаются числовые характеристики тождеств 
конечномерных неассоциативных алгебр. Основным результатом является построение
четырехмерной простой алгебры с единицей, имеющей дробную PI-экспоненту,
строго меньшую ее размерности.}

\vskip .2in
{\em Ключевые слова}: тождества, коразмерности, дробный экспоненциальный рост.
\vskip .2in

\section{Введение}

В работе изучаются тождества неассоциативных алгебр  над полем
$F$ нулевой характеристики. Все необходимые сведения по тождествам
в алгебрах можно найти в книгах \cite{b}, \cite{gzbook}, \cite{zsss}. 

С каждой алгеброй $A$ можно связать целочисленную последовательность $\{c_n(A)\}$, $n=1,2,\ldots$, которую называют последовательностью коразмерностей (все необходимые определения мы дадим ниже). Если данная последовательность растет экспоненциально, как, например, в конечномерном случае, то возникает вопрос о существовании предела $\lim_{n\to \infty}\sqrt[n]{c_n(A)}$. Существование и целочисленность такого предела были доказаны ранее для всех ассоциативных PI-алгебр \cite{gz98}, \cite{gz99}, для конечномерных алгебр Ли \cite{grz}, \cite{z2002}, для конечномерных йордановых и альтернативных алгебр \cite{gzTAMS}, \cite{gshz2011} и целого ряда других. В бесконечномерном случае этот предел, называемый PI-экспонентой, может быть как дробной, так и целой величиной даже в классе алгебр Ли (см. \cite{MRZ},\cite{zm}, \cite{vzm}), а в классе супералгебр Ли недавно обнаружены конечномерные 
алгебры с дробной экспонентой \cite{gz}.

В общем случае PI-экспонента конечномерной алгебры не превосходит ее размерности и 
может быть сколь угодно близкой к единице \cite{gmza} \cite{gmz}. Однако для двумерных алгебр PI-экспонента либо равна двум, либо единице (тогда $c_n(A)\le n+1$), либо равна 0 (в этом случае $c_n(A)=0$ для всех достаточно больших $n$) \cite{gmzPAMS}. Для трехмерных алгебр вопрос о целочисленности открыт, однако известно, что либо $c_n(A)\ge 2^{n}$ асимптотически, либо эта последовательность полиномиально ограничена. Однако для трехмерной алгебры с единицей PI-экспонента всегда существует и целочисленна \cite{AL}. Минимальная известная до недавних пор размерность, в которой появляется дробная экспонента, равна пяти \cite{gmzCA}, а в случае супералгебр Ли - семи  \cite{gz}. В тезисах \cite{mss} анонсирован результат
о существовании 4-х мерной коммутативной алгебры с дробной экспонентой.

Особое место занимает вопрос о существовании и значении PI-экспоненты простой конечномерной алгебры над алгебраически замкнутым полем. В большинстве изученных классов алгебр (ассоциативных \cite{r}, лиевских \cite{grz}, йордановых, альтернативных и некоторых других \cite{gzTAMS}, \cite{gshz2011}) эта величина равна размерности самой алгебры. Первые примеры простых алгебр, у которых PI-экспонента строго меньше размерности приведены в \cite{gz}. Минимальная размерность среди алгебр, предъявленных в \cite{gz}, равна 17. Но остался нерешенным вопрос об их целочисленности.

В данной работе построен пример четырехмерной простой алгебры с единицей с дробной 
PI-экспонентой (теорема \ref{t1}).

\section{Основные понятия и определения}

На протяжении всей статьи $F$ означает поле нулевой характеристики,
и все алгебры рассматриваются над $F$. При этом мы будем придерживаться 
соглашения, что в левонормированном произведении скобки опускаются, т.е.
$abc= (ab)c$.

Обозначим через $F\{X\}$
свободную неассоциативную алгебру над $F$ со счетным
множеством порождающих $X$. 
Напомним, что если $A$ - некоторая алгебра над $F$, а $f=f(x_1,\ldots,x_n)$ ---
неассоциативный полином из $F\{X\}$, то $f$ называют тождеством
$A$, если $f(a_1,\ldots,a_n)=0$ для любых $a_1,\ldots,a_n\in A$.
Совокупность всех тождеств алгебры $A$ образует идеал в $F\{X\}$,
который мы будем обозначать как $Id(A)$.

Обозначим через $P_n=P_n(x_1,\ldots,x_n)$ подпространтсво всех
полилинейных полиномов от $x_1,\ldots,x_n$ в $F\{X\}$. Тогда
 $P_n\cap Id(A)$ является подпространством всех 
полилинейных тождеств алгебры$A$ от $n$ переменных.

Напомним, что $n$-й коразмерностью тождеств алгебры $A$
(или просто $n$-й коразмерностью $A$) называется величина
$$
c_n(A)=\dim\frac{P_n}{P_n\cap Id(A)}.
$$

В широком ряде случаев последовательность $\{ c_n(A)\}$
растет не быстрее экспоненциальной функции от $n$, т.е.
$$
c_n(A) < a^n
$$
для некоторого $a>1$. Тогда последовательность коразмерностей
удовлетворяет условиям
$$
0\le \sqrt[n]{c_n(A)} \le a,
$$
и это позволяет дать следующее определение.

{\bf Определение.} Нижней PI-экспонентой алгебры $A$ называется нижний предел
$$
\underline{exp}(A)=\underline{\lim}_{n\to \infty}\sqrt[n]{c_n(A)}.
$$
Верхний предел
$$
\overline{exp}(A)=\overline{\lim_{n\to \infty}}\sqrt[n]{c_n(A)}
$$
называется верхней PI-экспонентой алгебры $A$. Если $\underline{exp}(A)$
и $\overline{exp}(A)$ равны, т.е. последовательность
 $\{\sqrt[n]{c_n(A)}\}$ имеет обычный предел, то он
называется PI-экспонентой алгебры $A$,
$$
exp(A)=\overline{exp}(A)=\underline{exp}(A).
$$

Пространство $P_n$ наделено естественным действием группы подстановок $S_n$:
$$
\sigma f(x_1,\ldots, x_n)=f(x_{\sigma(1)},\ldots, x_{\sigma(n)}),
$$
которое превращает $P_n$ в модуль над групповым кольцом $FS_n$, или,
для краткости, в $S_n$-модуль. При этом $P_n\cap Id(A)$ является
$S_n$-подмодулем в $P_n$. Напомним, что степень характера $\chi(M)$
некоторого $S_n$-модуля $M$ -- это размерность $M$, т.е.
$\deg\chi(M)=\dim M$. Кроме того, неприводимые характеры, т.е.
характеры неприводимых представлений, однозначно задаются разбиениями
$\lambda\vdash n$ числа $n$, где $\lambda=
(\lambda_1,\ldots,\lambda_k)$, $\lambda_1\ge\ldots\ge \lambda_k>0$ и
$\lambda_1+\cdots+\lambda_k=n$ (все необходимые сведения по теории
представлений симметрических групп можно найти в \cite{J}, а по
ее применению в теории тождеств -- в \cite{b}, \cite{gzbook}).

Любой $S_n$-модуль $M$ раскладывается в сумму неприводимых, и в
терминах характеров это можно записать следующим образом:
$$
\chi(M)=\sum_{\lambda\vdash n} m_\lambda\chi_\lambda,
$$
где $\chi_\lambda$ - характер неприводимого модуля,
соответствующего $\lambda$, а $m_\lambda$ - его кратность в
разложении $M$. Отсюда
$$
\deg\chi(M)=\sum_{\lambda\vdash n} m_\lambda\deg \chi_\lambda.
$$
В частности, для $S_n$-модуля $\frac{P_n}{P_n\cap Id(A)}$ получаем
\begin{equation}\label{eq1}
\chi_n(A)=\chi\left(\frac{P_n}{P_n\cap Id(A)}\right)=
\sum_{\lambda\vdash n} m_\lambda\chi_\lambda,
\end{equation}
и
\begin{equation}\label{eq2}
c_n(A)=\sum_{\lambda\vdash n} m_\lambda\deg\chi_\lambda.
\end{equation}
Разложение (\ref{eq1}) называется $n$-м кохарактером $A$, а формула
(\ref{eq2}) позволяет получать оценки величины $c_n(A)$. Будем
говорить, что характер $\chi(M)=\sum_{\lambda\vdash n}
m_\lambda\chi_\lambda$ модуля $M$ лежит в полосе ширины $d$, если
все разбиения $\lambda$ с ненулевыми кратностями $m_\lambda$ имеют
не более $d$ частей, т.е. $\lambda=(\lambda_1,\ldots,\lambda_k)$ и
$k\le d$. Для таких разбиений известно, что $\deg\chi_\lambda\le d^n$.

Если $A$ - конечномерная $F$-алгебра, $\dim A=d$, то известно
(см., например, \cite[Теорема 4.6.2]{gzbook} и \cite{gmz}), что ее кохарактер
$\chi_n(A)$ лежит в полосе ширины $d$ и что сумма
$$
l_n(A)=\sum_{\lambda\vdash n} m_\lambda,
$$
называемая $n$-й кодлиной, в (\ref{eq1}) ограничена полиномиальной функцией 
от $n$. Более точно, в теореме 1 из \cite{gmz} доказано, что
\begin{equation}\label{eq3}
l_n(A)\le d(n+1)^{d^2+d}.
\end{equation}

Это, в частности, означает, что при экспоненциальном росте коразмерностей
ключевую роль играют только максимальные размерности $\deg\chi_\lambda$.
Для получения их оценок удобно
использовать следующие функции $\Phi(\lambda)$, связанные с
разбиениями $\lambda\vdash n$. Пусть
$\lambda=(\lambda_1,\ldots,\lambda_k)\vdash n$, $k\le d,
\lambda_1+\ldots+\lambda_k=n$. Положим
$$
\Phi(\lambda)=\frac{1}{\left( \frac{\lambda_1}{n}  \right)^{ \frac{\lambda_1}{n}}
\ldots \left( \frac{\lambda_d}{n}  \right)^{ \frac{\lambda_d}{n}}}\quad .
$$
(Если $k$ строго меньше $d$, то соответствующие множители $0^0$ равны $1$).

В работе \cite{gz} отмечалась связь между значением функции $\Phi(\lambda)$
и степенью характера $\deg\chi(\lambda)$ (см. лемму 1). Пусть $\lambda=
(\lambda_1,\ldots, \lambda_k)\vdash n$, $n\ge 100$ и $k\le d$. Тогда
\begin{equation}\label{eq0}
\frac{\Phi(\lambda)^n}{n^{d^2+d}} \le \deg\chi_\lambda \le n \Phi(\lambda)^n.
\end{equation}

Пусть теперь $A$ - произвольная конечномерная алгебра над $F$. Рассмотрим ее
$n$-й кохарактер (\ref{eq1}) и обозначим через $\Phi^{(n)}_{max}$ максимальное значение $\Phi(\lambda)$ по всем $\lambda\vdash n$,  у которых $m_\lambda\ne 0$ в
(\ref{eq1}). Комбинируя соотношения (\ref{eq2}), (\ref{eq3}) и (\ref{eq0}), 
мы получаем следующее утверждение.

\begin{lemma} \label{l1}
Пусть $\dim A=d$. Тогда 
$$
\frac{1}{n^{d^2+d}}(\Phi^{(n)}_{max})^n\le c_n(A) \le (n+1)^{d^2+d+1} (\Phi^{(n)}_{max})^n
$$
для всех $n\ge 100d$.
\end{lemma}
\hfill $\Box$

Нам также потребуется следующее свойство функции $\Phi(\lambda)$.
Если $\lambda=(\lambda_1,\ldots, \lambda_q)$ и $\mu=(\mu_1,\ldots,\mu_q)$ - два разбиения одного и того же целого числа $n$, то диаграммой Юнга $D_\lambda$, отвечающей разбиению $\lambda$, называется таблица из $\lambda_1$ клеток в 1-й строке, $\lambda_2$ клеток во 2-й строке, и т.д. Аналогично строится и диаграмма $D_\mu$
по разбиению $\mu\vdash n$. Мы будем говорить, что диаграмма $D_\mu$ получена из диаграммы $D_\lambda$ выталкиванием вниз одной клетки, если существуют такие $1\le i<j\le q$, что $\mu_i=\lambda_i-1$, $\mu_j=\lambda_j+1$ и $\mu_p=\lambda_p$ для всех
остальных $1\le p\le q$. Если же $\lambda=(\lambda_1,\ldots, \lambda_q)\vdash n$, а $\mu=(\mu_1,\ldots,\mu_q,1)\vdash n$, то $D_\mu$ получена из  $D_\lambda$ выталкиванием вниз одной клетки, если одна из строк $D_\mu$ на одну
клетку короче, чем у $D_\lambda$, а все остальные (кроме последней) имеют такую же длину.

\begin{lemma} \label{lq}
Пусть $D_\mu$ получена из $D_\lambda$ выталкиванием вниз одной клетки. Тогда
$\Phi(\mu)\ge \Phi(\lambda)$.
\end{lemma}

{\em Доказательство.} Пусть $\lambda=(\lambda_1,\ldots, \lambda_q)$ и $\mu=(\mu_1,\ldots,\mu_{q'})$ - два разбиения $n$, где $q'=q$ или $q+1$. Тогда
\begin{equation}\label{eqn}
\Phi(\lambda)^n=\frac{n^n}{\lambda_1^{\lambda_1}\ldots \lambda_q^{\lambda_q}}
\end{equation}
Если $q=q'$, то в аналогичном выражении для $\Phi(\mu)^n$ знаменатель получается из знаменателя в (\ref{eqn}) заменой одного произведения вида $a^ab^b, a\ge b+2$, на
$(a-1)^{a-1} (b+1)^{b+1}$. В этом случае неравенство $\Phi(\mu)^n\ge \Phi(\lambda)^n$
следует из того, что функция $f(x)=x^x  (c-x)^{c-x}$ убывает на интервале 
$(c/2; 0)$. Если же $q'=q+1$, то в знаменателе (\ref{eqn}) мы заменяем
множитель $a^a$ на $(a-1)^{a-1}\cdot 1^1< a^a$, и снова получаем $\Phi(\mu)^n > \Phi(\lambda)^n$. Поскольку $\Phi(\lambda), \Phi(\mu)>0$, отсюда следует 
утверждение леммы.
\hfill $\Box$

\section{Четырехмерная алгебра и ее кохарактер}

Рассмотрим четырехмерное векторное пространство $W$ с базисом $\{e_{-1}, e_0,e_1,e_2\}$ и введем на нем умножение следующим образом:
\begin{itemize}
\item[(a)]
$e_ie_0=e_0e_i=e_i$ для всех $-1\le i\le 2$;
\item[(б)] 
при $i,j \ne 0$ положим $e_ie_j=0$, если $i>j$ или $i+j<-1$, или $i+j>2$;
\item[(в)]
$e_ie_j=e_{i+j}$ во всех остальных случаях.
\end{itemize}

Легко заметить, что $W$ - простая алгебра с единицей $e_0$ и что разложение
$$
W=W_{-1}\oplus W_0\oplus W_1\oplus W_2,
$$
где $W_i=<e_i>, i=-1,\ldots,2$, является ${\mathbb Z}$-градуировкой на $W$.

Рассмотрим кохарактер
\begin{equation}\label{e5}
\chi_n(W)=\sum_{\lambda\vdash n} m_\lambda \chi_\lambda
\end{equation}
алгебры $W$. Поскольку $\dim W=4$, ее кохарактер лежит в полосе ширины 4, т.е.
$\lambda_5=0$ для любого разбиения $\lambda \vdash n$, если $m_\lambda\ne 0$ в (\ref{e5}).

Напомним, что таблицей Юнга $T_\lambda$ называется диаграмма $D_\lambda$, в клетках которой расставлены числа $\{1,\ldots, n\}$. Любой неприводимый модуль над групповой алгеброй $R=FS_n$ изоморфен минимальному левому идеалу $Re_{T_\lambda}$, где $e_{T_\lambda}$ - элемент $R$, построенный по следующему правилу.

Назовем стабилизатором по строкам $R_{T_\lambda}$ подгруппу всех подстановок из $S_n$, перемещающих символы $1,\ldots,n$ только в пределах своих строк, а стабилизатором по столбцам $C_{T_\lambda}$ --- подгруппу, переставляющую символы внутри столбцов. Обозначим
$$
R(T_\lambda)=\sum_{\sigma\in R_{T_\lambda}} \sigma~,\quad
C(T_\lambda)=\sum_{\tau\in C_{T_\lambda}} ({\rm sgn}~\tau)\tau~,\quad
e_{T_\lambda}=R(T_\lambda)C(T_\lambda).
$$
Тогда $e_{T_\lambda}$ --- квазиидемпотент группового кольца $FS_n$, т.е.
$e_{T_\lambda}^2=\alpha e_{T_\lambda}, \alpha\ne 0$ и $FS_n e_{T_\lambda}$ ---
неприводимый $S_n$-модуль с характером $\chi_\lambda$. Более того, если $M$
--- некоторый $FS_n$-модуль и
$$
\chi(M)=\sum_{\lambda\vdash n} m_\lambda \chi_\lambda,
$$
то $m_\lambda\ne 0$ тогда и только тогда, когда $e_{T_\lambda}M\ne 0$.

Пусть теперь $f$ --- полилинейный многочлен степени $n$, порождающий неприводимый
$FS_n$-подмодуль $M$ в $P_n$ с характером $\chi_\lambda$. Можно считать, что 
$e_{T_\lambda}f\ne 0$ для некоторой таблицы Юнга $T_\lambda$. Тогда и многочлен
$g=C(T_\lambda)f$ тоже порождает $M$ как $R$-модуль.
Входящие в $g$ переменные разбиваются на $m$ непересекающихся  подмножеств 
$$
\{x_1,\ldots, x_n\}= X_1\cup\cdots\cup X_m,
$$
где $m$ - число столбцов диаграммы Юнга $D_\lambda$, а $X_j$ - набор переменных, номера
которых записаны в $j$-м столбце. При этом $g$ кососимметричен по каждому из наборов 
$X_1,\ldots, X_m$. Другими словами, любой неприводимый $FS_n$-подмодуль в $P_n$ с 
характером $\chi_\lambda$ порождается полилинейным многочленом, зависящим от $m$ 
непересекающихся кососимметричных наборов переменных мощности $\lambda_1',\ldots,\lambda_m'$
соответственно, где $\lambda'=(\lambda_1',\ldots,\lambda_m')$ - сопряженное к $\lambda$
разбиенини  (т.е. $\lambda_1',\ldots,\lambda_m'$ --- высоты столбцов $D_\lambda$).

Чтобы получить более сильные органичения (чем $\lambda_5=0$) на кохарактер алгебры $W$, 
введем для разбиений еще одну числовую характеристику. Пусть $\lambda= (\lambda_1,\ldots,
\lambda_k)\vdash n$. Запишем во все клетки 1-й строки число $-1$, второй - $0$, и т.д.,
то есть в клетках $k$-й строки стоит число $k-2$. Тогда весом диаграммы $D_\lambda$
назовем сумму всех чисел в $D_\lambda$. Эту же величину назовем весом $wt(\lambda)$ разбиения $\lambda$. Другими словами,
$$
wt(\lambda)=\sum_{i=1}^k (i-2)\lambda_i.
$$

\begin{lemma}\label{l3}
Если $m_\lambda\ne 0$ в разложении (\ref{e5}), то
\begin{itemize}
\item
 $\lambda_5=0$;
\item
$wt(\lambda)\le 2$, т.е. $\lambda_1-\lambda_3-2\lambda_4\ge -2$.
\end{itemize}
\end{lemma}

{\em Доказательство.} Равенство $\lambda_5=0$ уже показано ранее. Докажем второе утверждение 
леммы. Пусть $m_\lambda\ne 0$, т.е. существует неприводимый $FS_n$-подмодуль в $P_n$ с 
характером $\chi_\lambda$, не лежащий в идеале тождеств алгебры $W$. Как отмечено выше, это
означает, что существует полилинейный многочлен $g=g(x_1,\ldots,x_n)$, кососимметричный
по переменным каждого из столбцов $T_\lambda$, не равный нулю тождественно в $W$. Но тогда
существует подстановка $\varphi: X\to \{e_{-1},e_0,e_1,e_2\}$, при которой $\varphi(g)\ne 0$.
В силу кососимметричности $g$ вместо переменных одного столбца мы вынуждены подставлять
различные базисные элементы алгебры $W$. В частности, если $x_{i_1},\ldots, x_{i_4}$ - 
переменные одного столбца высоты 4, то суммарный вес элементов $\varphi(x_{i_1}),\ldots, 
\varphi(x_{i_4})$ равен $-1+0+1+2=2$ в ${\mathbb Z}$-градуировке $W$. Поэтому минимальный
возможный вес $\varphi(g)$ в ${\mathbb Z}$-градуировке равен
$$
wt(\lambda)=-\lambda_1+\lambda_3+2\lambda_4.
$$
Поскольку все компоненты $W_k$ при $k\ge 3$ равны нулю, это и означает, что 
$-\lambda_1+\lambda_3+2\lambda_4\le 2$, и лемма доказана.
\hfill $\Box$

Нам потребуется в дальнейшем одно достаточное условие, при котором кратность не равна
нулю. Рассмотрим разбиение $\lambda=(\lambda_1,\lambda_2,\lambda_3,\lambda_4)$  и перепишем 
его в виде
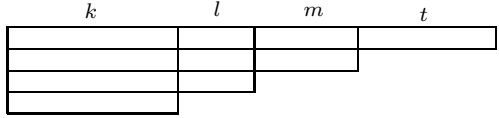
$$
\lambda=(k+l+m+t, k+l+m, k+l,k).
$$
Диаграмма этого разбиения выглядит следующим образом:

\vskip.1cm \hskip2cm
\setlength{\unitlength}{2565sp}%
\begingroup\makeatletter\ifx\SetFigFont\undefined%
\gdef\SetFigFont#1#2#3#4#5{%
  \reset@font\fontsize{#1}{#2pt}%
  \fontfamily{#3}\fontseries{#4}\fontshape{#5}%
  \selectfont}%
\fi\endgroup%
\begin{picture}(2865,2177)(464,-1286)
\thinlines {\color[rgb]{0,0,0}\put(476,364){\line( 1, 0){4725}}
\put(5201,364){\line( 0,-1){210}} \put(5201,148){\line(-1, 0){4727}}
\put(2871,-268){\line(-1, 0){2402}}
\put(3864,-66){\line(-1, 0){3388}} \put(2126,-474){\line(-1,
0){1650}} \put(476,-472){\line( 0, 1){838}}
}%
{\color[rgb]{0,0,0}\put(2126,364){\line( 0,-1){843}}
}%
{\color[rgb]{0,0,0}\put(3864,-71){\line( 0, 1){435}}
}%
{\color[rgb]{0,0,0}\put(2864,-271){\line( 0, 1){635}}
}%
\put(1226,451){\makebox(0,0)[lb]{\smash{{\SetFigFont{8}{9.6}{\rmdefault}{\mddefault}
{\updefault}{\color[rgb]{0,0,0}$k$}%
}}}}
\put(2471,476){\makebox(0,0)[lb]{\smash{{\SetFigFont{8}{9.6}{\rmdefault}{\mddefault}{\updefault}{\color[rgb]{0,0,0}$l$}%
}}}}
\put(3341,476){\makebox(0,0)[lb]{\smash{{\SetFigFont{8}{9.6}{\rmdefault}{\mddefault}{\updefault}{\color[rgb]{0,0,0}$m$}%
}}}}
\put(4464,426){\makebox(0,0)[lb]{\smash{{\SetFigFont{8}{9.6}{\rmdefault}{\mddefault}{\updefault}{\color[rgb]{0,0,0}$t$}%
}}}}
\end{picture}%

При этом $n=4k+3l+2m+t$, вес разбиения $\lambda$ равен $-m-t+2k$, а небходимое условие для неравенства $m_\lambda\ne 0$ принимает вид
$$
m+t\ge 2k-2.
$$

\begin{lemma}\label{l4}
Пусть $m+t\ge 2k$ и $m\le 2k$. Тогда $m_\lambda\ne 0$ в разложении (\ref{e5}).
\end{lemma}

{\em Доказательство.} Чтобы доказать лемму, достаточно построить полилинейный полином $f$
степени $n$, зависящий от $k$ кососимметричных наборов переменных
$\{x_1^{(i)},\ldots, x_4^{(i)}\}, 1\le i \le k$,  $l$ косомимметричных наборов
$\{y_1^{(i)},y_2^{(i)},y_3^{(i)}\}, 1\le i \le l$, $m$ косомимметричных наборов
$\{z_1^{(i)},z_2^{(i)}\}, 1\le i \le m$, и еще $t$ переменных  $u_1^{(i)}$, $1\le i\le t$,
такой, что после его симметризации по наборам $\{x_1^{i_1}, y_1^{i_2}, z_1^{i_3}, u_1^{i_4}\}$,
$\{x_2^{i_1}, y_2^{i_2}, z_2^{i_3}\}$, $\{x_3^{i_1}, y_3^{i_2}\}$ и по $\{x_1^{i_4}\}$, где
$1\le i_1 \le k, 1\le i_2 \le l, 1\le i_3 \le m, 1\le i_4 \le t$, получается полином, не равный нулю тождественно в $W$.

Для обозначения кососимметричных наборов удобно помечать входящие в него переменные
одним и тем же символом сверху. Например,
$$
\bar x_1 \bar x_2 \bar x_3  = \sum_{\sigma \in S_3} 
({\rm sgn}~\sigma) x_{\sigma(1)} x_{\sigma(2)} x_{\sigma(3)},
$$
$$
\bar a_1\widetilde b_1
 \bar a_2\widetilde b_2
 = a_1b_1a_2b_2 - a_1b_2a_2b_1 
- a_2b_1a_1b_2 + a_2b_2a_1b_1,
$$
а
$$
(\bar x \bar{\bar x}) (\bar y \bar{\bar y}) = 
(xx)(yy) - (yx)(xy) - (xy)(yx)+(yy)(xx).
$$

Мы будем использовать данное соглашение и для элементов алгебры $W$, а не только для переменных.
Например,
$$
\bar e_{-1} (\bar e_{1} \bar e_{2}) = e_{-1} (\bar e_{1} \bar e_{2}) - e_1(\bar e_{-1}\bar e_2) -
e_2(\bar e_{1} \bar e_{-1})=0-e_1(e_{-1}e_2)+e_2(e_{-1}e_1)=
$$
$$
=-e_1^2+e_2e_0=-e_2+e_2=0.
$$

Сначала построим выражение, кососимметричное по $e_{-1},e_0,e_1,e_2$  и содержащее два
дополнительных базисных элемента $e_{-1}$. Положим
$$
f_1=e_{-1}[\bar e_{-1}((\bar e_0e_{-1})(\bar e_1 \bar e_2))].
$$
Этот элемент является значением полинома

$$
x_{-1}[\bar x_{-1}((\bar x_0x_{-1})(\bar x_1 \bar x_2))]
$$
от переменных $x_{-1},x_0,x_1,x_2$ и равен 
$$
f_1=e_{-1}[\bar e_{-1}((e_0e_{-1})(\bar e_1 \bar e_2))] = e_{-1}[\bar e_{-1}(e_{-1}(\bar e_{1}\bar e_{2}))] =
$$
$$
= e_{-1}[ e_{-1}(e_{-1}(\bar e_{1}\bar e_{2}))-e_{1}(e_{-1}(\bar e_{-1}\bar e_{2})) 
-e_{2}(e_{-1}(\bar e_{1}\bar e_{-1}))]=
$$
$$ 
=e_{-1}[0-e_{1}(e_{-1}(e_{-1}e_2))+e_2(e_{-1}(e_{-1}e_1))]=
e_{-1}[-e_1(e_{-1}e_1)+e_2(e_{-1}e_0)]=
$$
$$
=e_{-1}[-e_1e_0+e_2e_{-1}=e_{-1}[-e_1+0]=-e_0.
$$

Обозначим
$$
g(x_{-1}, x_0,x_1,x_2,y,z)=y[\bar x_{-1}((\bar x_0z)(\bar x_1 \bar x_2))].
$$
Тогда левонормированная степень $f_1^k$ является значением симметризации многочлена
$$
g(x_{-1}^{(1)}, x_{0}^{(1)}, x_{1}^{(1)},x_{2}^{(1)},y_{1}^{(1)}
, z_{1}^{(1)})
\ldots
g(x_{-1}^{(k)}, x_{0}^{(k)}, x_{1}^{(k)},x_{2}^{(k)},y_{1}^{(k)}
, z_{1}^{(k)}),
$$
косомисметричного по $\{x_{-1}^{(i)}, x_{0}^{(i)}, x_{1}^{(i)},x_{2}^{(i)},\}$, $i=1,\ldots,k$. Симметризация же проведена по четырем наборам
$$
\{x_{-1}^{(1)},\ldots, x_{-1}^{(k)}, y_{1}^{(1)}
, z_{1}^{(1)}
,\ldots, y_{1}^{(k)}
, z_{1}^{(k)}
\},
\{x_{0}^{(1)},\ldots, x_{0}^{(k)}\}
, 
$$
$$
\{x_{1}^{(1)},\ldots, x_{1}^{(k)}\}
, \{x_2^{(1)},\ldots, x_2^{(k)}\}
.
$$
Аналогично, значением 3-альтерированного полинома является
$$
f_2=\bar e_{-1}  \bar e_{0} \bar e_{1}= (\bar e_{-1} \bar e_{0}) e_{1} +(\bar e_{0}
\bar e_{1}) e_{-1}+(\bar e_{1} \bar e_{-1}) e_{0} = 0+0-e_{-1}e_1e_0=-e_0.
$$
Положим также
$$
f_3=e_{-1}[ \bar e_{-1}((\bar e_{0}\widetilde e_{-1})(\bar e_{1}\bar e_2))]\widetilde e_0.
$$
Это выражение содержит кососимметричный набор $e_{-1}, e_0,e_1, e_2$, отмеченный чертой сверху, и
кососимметричный набор $e_{-1}, e_0$, помеченный волной. Поскольку $W_+e_{-1}=0$, где 
$W_+=W_1\oplus W_2$, то
$$
f_3=f_1e_0=f_1=-e_0.
$$

Наконец, пусть
$$
f_4=\bar{\bar e}_{-1}[[\bar e_{-1}((\bar e_0\widetilde e_{-1})(\bar e_1 \bar e_2))]\widetilde e_0
\bar{\bar e}_0].
$$
Здесь мы имеем один кососимметричный набор $\{e_{-1}, e_0,e_1, e_2\}$ и два кососимметричных набора $\{e_{-1}, e_0\}$. Обозначим
$$
a=[\bar e_{-1}((\bar e_0 \widetilde e_{-1})(\bar e_1 \bar e_2))]\widetilde e_0.
$$
Тогда степень $a$  в ${\mathbb Z}$-градуировке равна $1$, и поэтому из условия  $W_+e_{-1}=0$ и вычислений, проведенных выше для $f_1$ следует, что
$a=-e_1$ и
$$
f_4=e_{-1}(ae_0)-e_0(a e_{-1})=e_{-1}a=-e_0.
$$

Сначала рассмотрим частный случай $m=2k, t=0$. Положим
\begin{equation}\label{d1}
f=f(\underbrace{e_{-1},\ldots, e_{-1}}_\alpha, \underbrace{e_{0},\ldots, e_{0}}_\beta,
\underbrace{e_{1},\ldots, e_{1}}_\gamma, \underbrace{e_{2},\ldots, e_{2}}_\delta)=
(f_2^l)(f_4)^k.
\end{equation}
В записи $f$ элемент $e_{-1}$ встречается в $k$ кососимметричных наборах порядка 
$4$, в $2k=m$
 кососимметричных наборах порядка $2$  и в $l$ кососимметричных наборах порядка $3$. Неальтернируемые
$e_{-1}$ в $f$ отсутствуют. Общая степень $\alpha$ по $e_{-1}$ у $f$ равна $k+l+m$. Элемент
$e_0$ встречается тоже $k$ раз в $4$-альтернированных наборах, $m=2k$ раз ---  в $2$-альтернированных и
$l$ раз --- в $3$-альтернированных наборах, всего $k+l+m$ раз. Наконец, $e_1$  входит в $k$ 
кососиметричных наборов порядка $4$ и в $l$ кососимметричных наборов порядка $3$ а $e_2$ --- в $k$
кососимметричных наборов порядка $4$. Это и означает, что $f$ является значением полинома, порождающего
неприводимый модуль, соответствующий разбиению $\lambda=(k+l+m+t, k+l+m, k+l,k)$, где $m=2k, t=0$.
Следовательно, $m_\lambda\ne 0$ для этого разбиения, т.к. $f=(-e_0)^{k+l}=\pm e_0$. 

Рассмотрим более общий случай, когда $m=2k$ и $t>0$. Положим $n_0=n-t$ и для разбиения
$\lambda_0=(k+l+m, k+l+m, k+l,k)\vdash n_0$ построим многочлен $f$, указанный в $(\ref{d1})$. Тогда
$$
g=f(\underbrace{e_{-1}+e_0,\ldots, e_{-1}+e_0}_\alpha, e_0,\ldots,e_2)\underbrace{ (e_{-1}+e_0) \cdots (e_{-1}+e_0)}_t=
$$
$$
=f(e_{-1}+e_0,\ldots, e_{-1}+e_0, e_0,\ldots,e_2)+f',
$$
где $f'=tf(e_{-1},\ldots,e_{-1},e_0,\ldots,e_2)e_{-1}\in W_{-1}$
поскольку $(W_{-1}\oplus W_1\oplus W_2)e_{-1}=0$. Кроме того,
$$
f(e_{-1}+e_0,\ldots, e_{-1}+e_0, e_0,\ldots,e_2)= f(e_{-1},\ldots, e_{-1}, e_0,\ldots,e_2) + f'',
$$
где $f''\in W_1\oplus W_2$. Отсюда следует, что $g=\pm e_0+f'+f''\ne 0$, и $g$ является значением полинома,
порождающего неприводимый $FS_n$-модуль с характером $\chi_\lambda$. Следовательно, $m_\lambda\ne 0$ 
в (\ref{e5}).

Теперь пусть $m=2q < 2k$ и $t\ne 0$. Рассмотрим произведение
$$
(f_1)^{k-q}(f_2)^l(f_4)^q=f_0=f_0(\underbrace{e_{-1},\ldots, e_{-1}}_\alpha, \underbrace{e_{0},\ldots, 
e_{0}}_\beta,\underbrace{e_{1},\ldots, e_{1}}_\gamma, \underbrace{e_{2},\ldots, e_{2}}_\delta).
$$
Как и выше, $\delta=k$, и все $e_2$ входят в кососимметричные наборы порядка $4$, $\gamma=k+l,\beta=
k+l+m$, и все $e_1,e_0$ входят в кососимметричные наборы, причем $k$ из них для $e_1$ - порядка $4$,
а $l$ - порядка $3$. Для $e_0$ число наборов порядка $4,3$ и $2$ равно $k,l$ и $m$ соответственно.
Элемент $e_{-1}$ тоже встречается в кососимметричных наборах порядка $4,3$ и $2$, число которых равно
$k-q+q=k, l$ и $2q=m$ соответственно. Но кроме этого, $e_{-1}$ встречается еще $2(k-q)$ раз вне 
кососимметричных наборов за счет множителя $f_1^{k-q}$. В частности,
$$
\alpha=k+l+m+t_0,
$$
где $t_0=2(k-q)$, т.е. $m+t_0=2k$.

Раccмотрим выражение
$$
f_0'=f_0(\underbrace{e_{-1}+e_0,\ldots, e_{-1}+e_0}_\alpha, \underbrace{e_{0},\ldots, 
e_{0}}_{\beta=k+l+m}, \underbrace{e_{1},\ldots, e_{1}}_{\gamma=k+l}, 
\underbrace{e_{2},\ldots, e_{2}}_{\delta=k}),
$$
 а вместе с ним и выражение
$$
f=f_0'\underbrace{(e_{-1}+e_0)\ldots (e_{-1}+e_0)}_{t-t_0}.
$$

Как и в предыдущем случае,
$$
f=f_0+f'+f'',
$$
где $f_0=\pm e_0$, $f'\in W_{-1}$ а $f''\in W_{+}$, т.е. $f\ne 0$. Тогда $f$ является ненулевым значением 
полинома, отвечающего разбиению $\lambda=(k+l+m+t, k+l+m, k+l,k)$ с $m=2q< 2k$.
Значит и для таких разбиений $\lambda$ кратности в (\ref{e5}) ненулевые.

Наконец, при нечетном $m=2q+1<2k$, мы сначала берем произведение
$$
f_0=f_1^{k-q-1} f_2^l f_3 f_4^q,
$$
а затем заменяем в нем $e_{-1}$ на $e_{-1}+e_0$ и берем
$$
f=f_0(e_{-1}+e_0,\ldots,e_{-1}+e_0,e_0,\ldots, e_0,\ldots, e_2,\ldots e_2)
\underbrace{(e_{-1}+e_0)
\ldots (e_{-1}+e_0)}_{t-t_0},
$$
где $t_0=2(k-q)-1$. Тогда $e_{-1}'=e_{-1}+e_0$ в записи $f_0$ входит в $k$
косомсимметричных наборов из четырех элементов, в $l$ наборов из 3-х элементов, в $m=2q+1$
наборов из 2-х элементов
 и вне альтернированных наборов всртечается $t$ раз. Для 
$e_0,e_1,e_2$ выполняются те же условия, что и в предыдущем случае. Поскольку $f_0=\pm e_0$,
а $f=f_0+f'$, где $f'\in W_{-1}\oplus W_1\oplus W_2$, то $f\ne 0$. Другими словами, кратность $m_\lambda$ не равна
нулю и для $\lambda=(k+l+m+t, k+l+m, k+l,k)$ с нечетным $m<2k$, если $m+t\ge 2k$. Таким образом, для $k\ne 0$ лемма доказана.

Если же $k=0$, то и $m=0$, и разбиение $\lambda$ имеет вид $\lambda=(l+t,l,l)$.
Для него полилинейный многочлен строится аналогичным образом. Сначала возьмем
$f=f_2^l=f(e_{-1},\ldots,e_{-1}, e_0,\ldots,e_0,e_1,\ldots, e_1)$ степени $l$
по каждому из базисных элементов $e_{-1},e_0,e_1$. Затем положим
$$
f'=f(e_{-1}+e_0,\ldots,e_{-1}+e_0,e_0,\ldots, e_1)
\underbrace{(e_{-1}+e_0)\ldots (e_{-1}+e_0)}_t.
$$
Тогда, как и прежде, $f'$ не является тождеством $W$ и порождает $S_n$-модуль
с характером $\chi_\lambda$.

\hfill $\Box$

\section{Оценки PI-экспоненты}

Поскольку алгебра $W$ проста, то ее PI-экспонента существует (см. \cite[теорема 3]{gz}). Для 
ее оценки, в частности, для доказательства ее нецелочисленности введем следующую величину.

Пусть $\lambda=(\lambda_1,\ldots,\lambda_4)\vdash n$ и $m_\lambda$ --- кратность $\lambda$ в кохарактере (\ref{e5}). Положим
$$
a_n=\max\{\Phi(\lambda)|m_\lambda\ne 0\}.
$$
Тогда в силу леммы \ref{l1}
\begin{equation}\label{eq6}
\overline{exp(A)}= \overline{\lim}_{n\to \infty} a_n,\quad 
\underline{exp(A)}= \underline{\lim}_{n\to \infty} a_n
\end{equation}
для любой четырехмерной алгебры. 

Рассмотрим теперь последовательность $\{a_n\}$ для алгебры $W$. Тогда в силу (\ref{eq6}) и \cite[теорема 3]{gz} последовательность $\{a_n\}$ имеет предел при $n\to \infty$.  Нам потребуется еще одно свойство этой последовательности.

\begin{lemma}\label{l5}
Пусть $a_n=\Phi(\lambda^{(0)})$. Тогда можно выбрать разбиение $\lambda^{(0)}$ таким, что $wt(\lambda^{(0)})\ge 0$ для достаточно больших $n$.
\end{lemma}

{\em Доказательство.} Пусть $\lambda^{(0)}$ - одна из точек максимума $\Phi(\lambda)$, 
определяющих $a_n$. Можно считать, что для всех таких $\lambda\vdash n$ с $\Phi(\lambda)=a_n$ 
это разбиение максимального веса. Как и раньше, запишем $\lambda^{(0)}$ в виде
$\lambda^{(0)}=(k+l+m+t,k+l+m,k+l,k)$. Предположим, что $wt(\lambda^{(0)})=-m-t+2k<0$, т.е.
$m+t>2k$.

Заметим сначала, что $k\ne 0$ у $\lambda^{(0)}$. Действительно, нетрудно заметить, что 
$\Phi(\lambda^{(0)})\le 3$ при $k=0$. В то же время для разбиения $\lambda=(3p,p,p,p)$
мы имеем
$$
\Phi(\lambda)=((\frac{1}{2})^{\frac{1}{2}}(\frac{1}{6})^{\frac{3}{6}})^{-1}
=\sqrt{12}>3,4.
$$
Поскольку последовательность $\{a_n\}$ сходится, а для любого $n$ найдется такое $p$, что
$|n-6p|\le 5$, то $a_n>3$ при всех достаточно больших $n$, и $k\ne 0$.

Заметим теперь, что если в диаграмме $D_\lambda$ перенести одну клетку из 2-й строки в 3-ю,
то получим $D_\mu$, где $\mu=(k'+l'+m'+t',k'+l'+m', k'+l',k')$ с $k'=k, t'=t+1$ и $m'=m-2$.
Поэтому либо $m'-2k'>0$, и тогда $m'+t'-2k'>0$, либо $m'-2k'=0$ или $-1$. В последних двух 
случаях выполняетcя неравенство $m'+t'-2k'\ge 0$. Таким образом, при $m\ge 2$ выталкивая вниз
одну или несколько клеток, мы переходим к разбиению $\mu$ большего веса, для которого
$\Phi(\mu)\ge \Phi(\lambda^{(0)})$ по лемме \ref{lq}. При этом $\mu$ удовлетворяет условиям леммы \ref{l4}, и
поэтому $m_\mu\ne 0$ в (\ref{e5}). В силу максимальности веса $\lambda^{(0)}$  мы получаем
условие $m\le 1$. 

Но тогда $t\ge 2k\ge 2$, и перебрасывая одну клетку у $D_{\lambda^{(0)}}$ из 1-й строки во 2-ю
мы получим диаграмму $D_\mu$, для которой $\Phi(\mu)\ge \Phi(\lambda^{(0)})$ и
$wt(\mu)> wt(\lambda^{(0)})$. При этом $\mu$ тоже удовлетворяет условиям леммы \ref{l4} и
$m_\mu\ne 0$. Отсюда следует, что $m+t-2k \le 0$ у $\lambda^{(0)}$, и лемма доказана.
\hfill $\Box$.

Ранее мы отмечали, что если диаграмма $D_\mu$ получена из диаграммы $D_\lambda$ выталкиванием 
вниз одной клетки, то $\Phi(\mu)\ge \Phi(\lambda)$. Теперь мы оценим насколько велико
это отклонение.

\begin{lemma}\label{l7}
Пусть $\lambda=(\lambda_1,\ldots, \lambda_q)$ и $\mu=(\mu_1,\ldots,\mu_{q'})$ - два разбиения числа $n$, где $q'=q$ или $q+1$ и $D_\mu$ получена из  $D_\lambda$ выталкиванием вниз одной клетки. Тогда
$$
\Phi(\lambda)\ge\frac{1}{n^{\frac{q^2+3q+4}{n}}}\Phi(\mu)
$$
\end{lemma}
{\em Доказательство.} Процедуру выталкивания клетки можно осуществить в два шага. Сначала
вырезаем клетку из $D_\lambda$ и получаем $D_{\lambda'}$, где $\lambda'\vdash n-1$,
а затем, подклеивая одну клетку к  $D_{\lambda'}$, получаем $D_\mu$. Тогда по лемме
6.2.4 из \cite{gzbook}
$$
\deg\chi_{\lambda'} \le \deg\chi_\lambda\le n \deg\chi_{\lambda'},\quad
\deg\chi_{\lambda'} \le \deg\chi_\mu\le n \deg\chi_{\lambda'},
$$
откуда легко вытекает, что
\begin{equation}\label{eq7}
\deg\chi_\lambda\ge \frac{1}{n} \deg\chi_\mu.
\end{equation}
Тогда из (\ref{eq7}) и формулы (\ref{eq0}) получаем
$$
\Phi(\lambda)^n \ge \frac{1}{n} \deg\chi_\lambda\ge \frac{1}{n^2} \deg\chi_\mu 
\ge \frac{1}{n^{(q+1)^2+q+1+2}}\Phi(\mu)^n,
$$
откуда и следует утверждение леммы.
\hfill $\Box$

Докажем еще одно соотношение, связывающее значения функции $\Phi$ для разных
разбиений.

\begin{lemma}\label{l7a}
Пусть диаграмма Юнга разбиения  $\lambda=(\lambda_1,\ldots, \lambda_d)\vdash (n-1)$
получена из диаграммы $D_\mu$ вычеркиванием одной клетки. Тогда
$$
\Phi(\lambda)\le n^\frac{d^2+d+2}{n}\Phi(\mu)
$$
при $n\ge d$.
\end{lemma}
{\em Доказательство.} 
По формуле (\ref{eq0})
$$
\Phi(\lambda)^{n-1}\le (n-1)^{d^2+d}\deg\chi_\lambda\le n^{d^2+d}\deg\chi_\lambda,
$$
а
$$
\deg\chi_\mu\le n\Phi(\mu)^n.
$$
С другой стороны, $\deg\chi_\lambda\le \deg\chi_\mu$ согласно \cite[лемма 6.2.4]{gzbook}. Поскольку максимальное значение $\Phi(\lambda)$ равно $d$, то отсюда получаем
$$
\Phi(\lambda)\le n^\frac{d^2+d+2}{n}\Phi(\mu).
$$
\hfill $\Box$

Определим еще одну последовательность, связанную с $W$. Пусть $n\ge 6$. Положим
$$
b_n=\max\{\Phi(\lambda)|\lambda=(\lambda_1,\ldots,\lambda_4)\vdash n, m_\lambda\ne 0,
\lambda_1-\lambda_3=2\lambda_4 \},
$$
если $n$ допускает разбиение $\lambda$ с $\lambda_1-\lambda_3=2\lambda_4$  и с
$m_\lambda\ne 0$ в (\ref{e5}) В противном случае положим $b_n=\min\{b_{n-1}, a_n\}$. Заметим, что согласно лемме \ref{l4} при
$n=6k$ разбиение $\lambda=(3k,k,k,k)$ удовлетворяет необходимым условиям.

\begin{lemma}\label{l9}
$$
\lim_{n\to \infty} b_n = \lim_{n\to\infty}a_n =exp(W).
$$
\end{lemma}
{\em Доказательство.} Согласно \cite[теорема 3]{gz} PI-экспонента любой конечномерной
простой алгебры существует, и следовательно существует предел
$\lim_{n\to\infty}a_n =exp(W)$ как вытекает из (\ref{eq6}). Теперь чтобы доказать лемму, достаточно найти такую функцию $\psi=\psi(n)$, что $\lim_{n\to\infty}\psi(n)=1$ и
\begin{equation}\label{eq8a}
\psi(n)a_n \le b_n\le a_n
\end{equation}
для всех достаточно больших $n$.

Зафиксируем $n$ и возьмем такое разбиение $\lambda\vdash n$, что $\Phi(\lambda)=a_n$.
По лемме \ref{l5} можно выбрать $\lambda$ так, что $wt(\lambda)\ge 0$. Но тогда из леммы \ref{l3} вытекает, что $wt(\lambda)$  равняется 0,1 или 2.

Если $wt(\lambda)=0$, то $b_n=a_n$. Пусть $wt(\lambda)=1$. Запишем $\lambda$ как
$\lambda=(k+l+m+t, k+l+m, k+l, k)$. Тогда $m+t= 2k-1$. Если $m\ne 0$, то можно перенести одну клетку из 2-й строки в 1-ю в диаграмме $D_\lambda$ и получить диаграмму $D_\mu$, $\mu=(k+l+m'+t',k+l+m', k+l, k)$, у которой
$m'=m-1, t'=t+2$, а $wt(\mu)=0$. Тогда по лемме \ref{l4}
$\mu$ имеет ненулевую кратность в $\chi_n(W)$. Как отмечалось в доказательстве леммы 
\ref{l5}, разбиение $\lambda$ имеет ненулевую компоненту $k$.
 Поэтому с учетом леммы \ref{l7}, мы имеем 
\begin{equation}\label{eq9}
b_n\ge \Phi(\mu)\ge \frac{\Phi(\lambda)}{n^\frac{32}{n}}=
\frac{a_n}{n^\frac{32}{n}}.
\end{equation}

Если $m=0$, но $l>0$ и $t>0$, то разбиение $\mu$ с нулевым весом можно получить переносом одной клетки из 3-й строки во 2-ю у $D_\lambda$, и мы снова получаем неравенство (\ref{eq9})  для $b_n$. Заметим, что случай $m=0$, $l>0$, $t=0$ невозможен, т.к. $m+t=2k-1$.

Единственное разбиение $\lambda$ с $wt(\lambda)=1$, которое не допускает указанных переносов, это разбиение $(3k-1,k,k,k)$. Но тогда
\begin{equation}\label{eq10}
\Phi(\lambda)\le n^\frac{22}{n}\Phi(\mu)
\end{equation}
по лемме \ref{l7a}, где $\mu=(3k,k,k,k)$. Поскольку $\Phi(\mu)=\sqrt{12}<3,48$,
то мы получаем
$$
\Phi(\lambda)< n^\frac{22}{n}\cdot 3,48.
$$

Заметим теперь, что разбиение вида $\rho=(3q,3q,q,q)$ удовлетворяет условиям леммы \ref{l4} и 
$$
\Phi(\rho)=\frac{8}{\sqrt[4]{27}}> 3,5,
$$
поэтому $\Phi(\lambda)$ не может удовлетворять неравенству (\ref{eq10}) при 
достаточно больших $n$, т.е. $\lambda\ne (3k-1,k,k,k)$, и неравенство (\ref{eq9})
выполнено, если $wt(\lambda)=1$.

Пусть теперь $wt(\lambda)=2$. Тогда будем поднимать два раза по одной клетке на
одну строку выше в диаграмме $D_\lambda$. Этого нельзя сделать только если
$\lambda=(3k-2,k,k,k)$, либо $\lambda=(q,q,q,1)$, либо после переноса одной клетки вверх мы получим разбиение $\mu=(3k-1,k,k,k)$. Первый и третий случаи исключаются также, как и при $wt(\lambda)=1$, поскольку такие разбиения не могут дать максимальне значение $\Phi(\lambda)$, а второй невозможен, поскольку $\deg\chi_\lambda\le n\deg
\chi_\mu$, если $\lambda=(q,q,q,1)\vdash n$ и $\mu=(q,q,q)\vdash (n-1)$, а
$\Phi(\mu)=3$.

Во всех остальных случаях двукратное применение леммы \ref{l7} дает нам соотношение
\begin{equation}\label{eq11}
b_n\ge \frac{a_n}{n^\frac{64}{n}}.
\end{equation}
Из (\ref{eq9}) и (\ref{eq11}) получаем требуемое условие (\ref{eq8a}), и лемма доказана.

\hfill $\Box$

Чтобы сформулировать и доказать основные результаты статьи, расширим область определения функции $\Phi$. Для любых $0\le x_1,\ldots, x_4\le 1$ положим
\begin{equation}\label{eq12}
\Phi(x_1,\ldots, x_4)=\frac{1}{x_1^{x_1}\ldots x_4^{x_4}}
\end{equation}
и рассмотрим внутри области определения $\Phi$ замкнутое подмножество $T$, заданное
условиями
\begin{equation} \label{eq13}
\left\{
               \begin{array}{ll}
        x_1\ge x_2\ge x_3\ge x_4     \\
x_1+x_2+x_3+x_4=1\\
x_1-x_3=2x_4.
               \end{array}
             \right.
\end{equation}

\begin{theorem}\label{t1}
PI-экспонента алгебры $W$ существует и равна
\begin{equation} \label{eq14}
exp(W)=\max\{\Phi(x_1,\ldots, x_4)\vert (x_1,\ldots, x_4)\in T\}.
\end{equation}
В частности, $exp(W)\approx3,610718614$.
\end{theorem}

{\em Доказательство.}Существование экспоненты уже отмечалось ранее, поскольку $W$ проста. Кроме того, 
$$
exp(W)=b=\lim_{n\to\infty}b_n
$$
по лемме \ref{l9}. Осталось показать, что $b=M$, где
$$
M=\max\{\Phi(x_1,\ldots, x_4)\vert (x_1,\ldots, x_4)\in T\}.
$$

Пусть $Z=(z_1,\ldots, z_4)$ --- точка максимума $\Phi$ на $T$. Ясно, что можно найти
точку $A=(a_1,\ldots, a_4)\in T$ с рациональными коэффициентами, сколь угодно близкую к $Z$. Обозначим через $m$ общий знаменатель рациональных чисел $a_1,\ldots, a_4$. Тогда $\lambda_1=a_1m,\ldots, \lambda_4= a_4m$ --- неотрицательные целые числа, 
причем $\lambda_1\ge\ldots \ge \lambda_4$. Другими словами, $\lambda=(\lambda_1, \ldots, \lambda_4)$ --- разбиение числа $m$, удовлетворяющее условию 
$\lambda_1-\lambda_3=2\lambda_4$. Более того, для любого $t=1,2,\ldots$ разбиение
$t\lambda=(t\lambda_1,\ldots, t\lambda_4)$ числа $n_t=tm$ удовлетворяет тому же условию. Отсюда следует, что
\begin{equation} \label{eq15}
b_{n_t}\ge \Phi(t\lambda)=\Phi(\lambda).
\end{equation}
Поскольку последовательность $\{b_i\}$ сходится, а $\Phi(\lambda)$ в (\ref{eq15})
можно сделать сколь угодно близкой к $M$, то $b\ge M$. Так как обратное условие очевидно, то мы доказали соотношение (\ref{eq15}).

Для полного завершения доказательства осталось обосновать приближенную оценку 
для $exp(W)$. В работе \cite{zm} был построен пример бесконечномерной алгебры Ли $L$, у которой $3,1< \underline{exp}(L)\le \overline{exp}(L)<3,9$. Недавно \cite{vzm}
было доказано существование обычной PI-экспоненты $L$, т.е. равенство
$\underline{exp}(L)= \overline{exp}(L)$. Более того, оказалось что
$$
exp(L)= \max\{\Phi(x_1,\ldots, x_4)\vert (x_1,\ldots, x_4)\in T\},
$$
где $\Phi$ --- функция, определенная в (\ref{eq12}), а область $T$ задана соотношениями (\ref{eq13}). В той же работе и было показано, что
$$
M=\Phi(\beta_1,\ldots,\beta_4),
$$
где $\beta_4$ ---  положительный корень уравнения $16t^3-24t^2+11t-1=0$,
$\beta_4\approx 0,276953179$, а
$$
\beta_3=2\beta_4-4\beta_4^2,\quad \beta_2=\frac{\beta_3^2}{\beta_4},
\quad \beta_1=\frac{\beta_3^3}{\beta_4^2}.
$$
Отсюда $exp(W)=exp(L)\approx 3,610718614$, и теорема доказана.
\hfill $\Box$

\begin{cor}
Существуют конечномерные простые алгебры с 1 с дробной экспонентой, строго
меньшей их размерности
\end{cor}
\hfill $\Box$

\begin{cor}
Наименьшая размерность алгебр с единицей, в которой встречаются дробные 
PI-экспоненты, равна 4.
\end{cor}
{\em Доказательство.}
Теорема \ref{t1} показывает, что в размерности 4 примеры существуют. То, что они отсутствуют среди двумерных и трехмерных алгебр, следует из результатов \cite{AL}.
\hfill $\Box$

\end{document}